\documentclass[journal,11pt]{article}

\usepackage{amsmath,amssymb,amsthm,amsfonts,mathrsfs, mathtools, bm, bbm, dsfont}
\usepackage{graphicx}

\usepackage{subfigure}
\usepackage{ifthen}
\usepackage{multirow, multicol}
\usepackage{float}
\usepackage{subeqnarray, array}
\usepackage{fullpage}
\usepackage[font=small]{caption}
\usepackage{xcolor}
 \graphicspath{{Figures/}}
\usepackage{subfigure}
\usepackage[hyphens]{url}
\usepackage{enumitem}
\usepackage{hyperref}
\usepackage{cite}
\usepackage{mdframed}
\usepackage{atbegshi}
\AtBeginDocument{\AtBeginShipoutNext{\AtBeginShipoutDiscard}}

\newcommand{\beq}{\begin{equation}}
\newcommand{\eeq}{\end{equation}}
\newcommand{\beqa}{\begin{eqnarray}}
\newcommand{\eeqa}{\end{eqnarray}}
\newcommand{\beqan}{\begin{eqnarray*}}
\newcommand{\eeqan}{\end{eqnarray*}}

\renewcommand{\(}{\left (}
\renewcommand{\)}{\right )}
\renewcommand{\[}{\left [}
\renewcommand{\]}{\right]}

\newcommand{\E}{\mathbb{E} }

\newcommand{\Nset}{\mathbb{N}}

\newcommand{\Acal}{{\cal A}}

\newcommand{\Ccal}{{\cal C}}

\newcommand{\Ncal}{{\cal N}}

\newcommand{\Pcal}{{\cal P}}

\newcommand{\Xcal}{{\cal X}}

\newcommand{\Gfrak}{\mathfrak{G}}

\newcommand{\argmax}{\mathop{\rm argmax}}

\renewcommand{\v}[1]{{\bm{#1}}}

\newcommand\T{{\mathpalette\raiseT\intercal}}
\newcommand\raiseT[2]{\raisebox{0.25ex}{$#1#2$}
}

\newcounter{l1}
\newcounter{l2}
\newcounter{l3}
\newcommand{\bdotlist}{\begin{list}{$\bullet$}{}}
\newcommand{\bboxlist}{\begin{list}{$\Box$}{}}
\newcommand{\bbboxlist}{\begin{list}{\raisebox{.005in}{{\tiny
$\blacksquare$ \ \ }}}{}}
\newcommand{\bdashlist}{\begin{list}{$-$}{} }
\newcommand{\blist}{\begin{list}{}{} }
\newcommand{\barablist}{\begin{list}{\arabic{l1}}{\usecounter{l1}}}
\newcommand{\balphlist}{\begin{list}{(\alph{l2})}{\usecounter{l2}}}
\newcommand{\bAlphlist}{\begin{list}{\Alph{l2}.}{\usecounter{l2}}}
\newcommand{\bdiamlist}{\begin{list}{$\diamond$}{}}
\newcommand{\bromalist}{\begin{list}{(\roman{l3})}{\usecounter{l3}}}

\newtheorem{theorem}{Theorem}

\newtheorem{proposition}{Proposition}

\newtheorem{definition}{Definition}
\newtheorem{assumption}{Assumption}

\renewcommand{\bar}{\overline}

\begin{document}

\title{Quantifying Market Efficiency Impacts of Aggregated Distributed Energy Resources}
{\thanks{K. Alshehri is with the Systems Engineering Department, King Fahd University of Petroleum and Minerals, Dhahran, 31261 Saudi Arabia. Email: kalshehri@kfupm.edu.sa. M. Ndrio, S. Bose, and T. Ba\c{s}ar are with the Department of Electrical and Computer Engineering, University of Illinois at Urbana-Champaign, Urbana, IL 61801 USA. Email: \{ndrio2, boses, basar1\}@illinois.edu.}}

{\author{Khaled Alshehri, Mariola Ndrio, Subhonmesh Bose, Tamer Ba\c sar\\
{{\it{\color{red}To appear in IEEE Transactions on Power Systems}}} }
\date{} 

\maketitle


\begin{abstract}
We focus on the aggregation of distributed energy resources (DERs) through a profit-maximizing intermediary that enables participation of DERs in wholesale electricity markets. Particularly, we study the market efficiency brought in by the large-scale deployment of DERs and explore to what extent such benefits are offset by the profit-maximizing nature of the aggregator. We deploy a game-theoretic framework to study the strategic interactions between an agreggator and DER owners. The proposed model takes into account the stochastic nature of the DER supply. We explicitly characterize the equilibrium of the game and provide illustrative examples to quantify the efficiency loss due to the strategic incentives of the aggregator. Our numerical experiments illustrate the impact of uncertainty and amount of DER integration on the overall market efficiency. 

\end{abstract}

\section{Introduction}

Widespread adoption of distributed energy resources (DERs) coupled with advances in communication and information technology, are pushing electricity markets to a more decentralized consumer-centric model. DERs typically include rooftop solar, small-scale wind turbines, electric vehicles as well as demand-response schedules. More generally, a DER is ``any resource on the distribution system that produces electricity and is not otherwise included in the formal NERC definition of the Bulk Electric System (BES)"\cite{derdefinition}. The low-voltage side of the grid, traditionally comprising mostly of passive small-scale consumers, is rapidly transforming into an active component of the grid where prosumers respond to price signals for managing their consumption and production of energy \cite{Morstyn}.

Transmission system operators lack visibility into the low and medium voltage distribution grid where DERs are connected. Furthermore, DERs have relatively small capacities that together with the high costs and complexities involved in their integration, render it impractical for such resources to directly offer their services in wholesale electricity markets. 
Despite significant research on effective means to harness such resources, a unifying framework for integration and compensation of DERs remains under debate.\footnote{FERC has issued Order 841 on energy storage and a Notice of Proposed Rulemaking on DERs, but is yet to define a binding framework for general DER participation in US wholesale markets.} See \cite{NYISO, DERbenefits,IEEEDER,NextGrid}  for insightful discussions.

One line of work promotes the implementation of distribution electricity markets operated by an independent distribution system operator (DSO), that acts as a market manager and dispatcher of DERs \cite{Lian,Ntakou,Manshadi}. In this model, the DSO is responsible for collecting the offers and bids from market participants and determine the appropriate prices to compensate DER asset owner-operators \cite{Sotkiewicz,Terra,Huang}. Another approach advocates fully distributed market structures, where prosumers trade DER services with each other as members of a coordinated and purely transactive community \cite{Moret,Rahimi}. 

In this paper, we focus on a third DER participation model through an aggregator $\Acal$ -- ``a company that acts as an intermediary between electricity end-users and DER owners, and the power system participants who wish to serve these end-users or exploit the services provided by these DERs", according to \cite{Burger}. California Independent System Operator (CAISO) allows such aggregators with capacities north of 0.5MW to participate in its wholesale markets for energy and ancillary services. See \cite{CAISODER} for an analysis of CAISO's model. Our focus is on efficiency impacts of a profit-motivated retail aggregator -- a topic that has not yet received much attention.

\begin{figure}
	\centering
	\includegraphics[width= .5\linewidth]{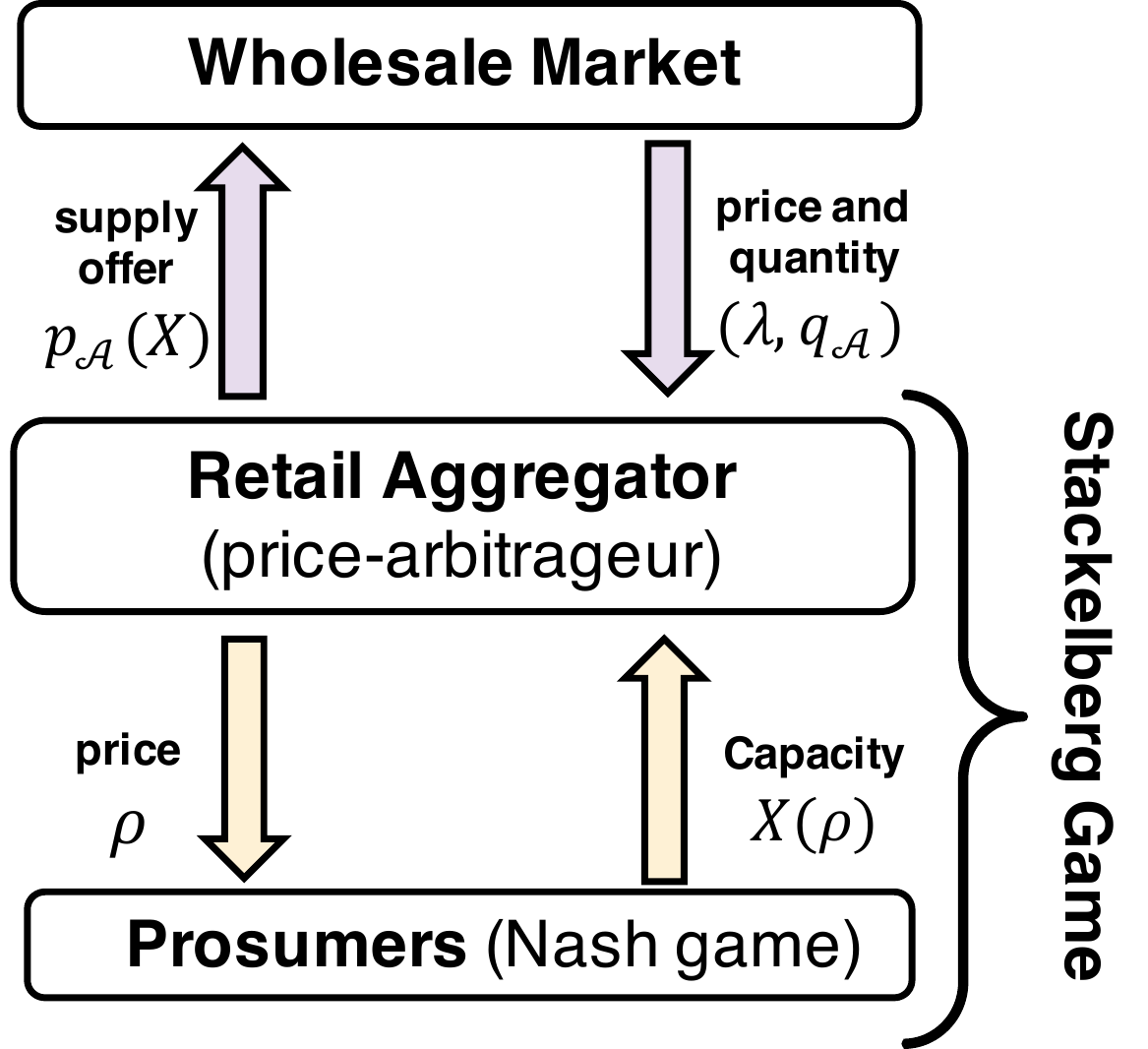}
	\caption{Interactions between prosumers, the DER aggregator, and the wholesale market.}
	\label{sketch_EW}
\end{figure}

We model the interaction between an aggregator $\Acal$ and prosumers with DERs as a Stackelberg game in Section \ref{sec:game}. $\Acal$ procures DER capacities from prosumers upon offering them a uniform price, and sells the aggregated capacity at the wholesale market price, profiting from arbitrage. DER capacities are uncertain. $\Acal$ faces penalties for defaulting on its promised offer to the wholesale market that she allocates among the prosumers. We analyze the resulting interaction, depicted in Figure \ref{sketch_EW} via game theory and characterize its equilibria that explicitly models uncertainties in DER supply. In Section \ref{sec:efficiency}, we introduce a new metric we call \textit{Price of Aggregation} ({\sf PoAg}) that compares power procurement costs in the wholesale market from two different DER participation models. In the first model, 
prosumers participate through $\Acal$. In the second one, they participate directly and offer their capacities in the wholesale market.  The second  model serves as the ideal yet impractical benchmark for efficient DER participation. {\sf PoAg} computes the efficiency loss due to the strategic nature of the aggregator. 
We analyze the game between $\Acal$ and the prosumers with various uncertainty models in DER capacities in Section \ref{sec:game2} and leverage these insights to study price of aggregation for illustrative example markets in Section \ref{sec:examples}. The results demonstrate how uncertainty and the amount of DER integration affect the {\sf PoAg} metric. We conclude the  paper in Section \ref{sec:conclusion} with remarks and future research directions. Proofs of all mathematical results are provided in the Appendix.

\section{The Game between the Prosumers and the DER Aggregator}\label{sec:game}

\newcommand{\LDA}{{\lambda_{\text{DA}}}}
\newcommand{\LRT}{{\lambda_{\text{RT}}}}

Consider a retail aggregator $\Acal$ who procures energy from a collection of prosumers $\Nset:=\{1,\ldots,N\}$ with DERs and offers the aggregate supply into the wholesale electricity market. $\Acal$ does not own any generation or consumption asset. She purely acts as an intermediary. To procure DER supply, she announces a uniform price $\rho$ at which she aims to buy energy from the prosumers. The latter respond by choosing how much energy each of the prosumers wishes to sell from their DERs such as rooftop photovoltaic panels, plug-in electric vehicles, wall-mounted batteries, thermostatically controlled loads, etc. DER supply is uncertain, owing to random variations in temperature, solar insolation, electric vehicle usage, etc. As a result, the realized aggregate supply from all prosumers may fall short of $\Acal$'s promised offer in the day-ahead (DA) market.\footnote{For ease of exposition, in this paper, we do not consider spatial variations in prices. However, our conclusions remain applicable to such considerations.} In that event, $\Acal$ faces a penalty for defaulting on its promise and allocates this penalty to the prosumers. In this section, we mathematize this interaction between $\Acal$ and the prosumers in $\Nset$ as a Stackelberg game (see Figure \ref{sketch_EW}). \footnote{For a comprehensive discussion on Stackelberg games, see to \cite{basar}.} Later in this paper, we explore how the outcomes of this game impact wholesale market efficiency.

Given a wholesale day-ahead market price $\LDA$,  $\Acal$ sets $\rho$, the price to procure energy from the prosumers. Then, prosumer $i \in \Nset$ responds by offering to sell $x_i$ to $\Acal$, who then  offers the aggregate procured DER supply\footnote{We do not consider economic witholding or strategic capacity reporting.}
$$ X := \v{1}^\T (x_1, \ldots, x_N)^\T := \v{1}^\T \v{x}$$ 
to the DA market at price $\LDA$. Here, $\v{1}$ is a vector of all ones of appropriate size. We assume that $\Acal$ is a price-taker in the wholesale market. That is, she believes that her wholesale offer will not influence the wholesale market prices. This assumption is natural as aggregators today typically do not command enough DER supply to exercise significant market power. In order to compute the DA offer, let $\Acal$ believe that its entire offer will be cleared in the DA market. Thus, $\Acal$ hopes to earn $(\LDA  - \rho) X$ from DA transactions. The revenue from DA sales stems from price arbitrage. $\Acal$ buys $X$ from prosumers at price $\rho$ and sells it to the wholesale market at price $\LDA$, bagging the difference. Setting a higher $\rho$ reduces the price difference from the DA market price, but generally incentivizes the prosumers to sell more of their DER supply, in turn, increasing the energy that $\Acal$ can offer in the DA market.

Assume that prosumers in $\Nset$ are homogenous, each of whom can supply power from a collection of DERs with installed capacity $\bar{C}$. Let $C_i\in[0,\bar{C}]$ denote the sum-total of capacities from all DERs with prosumer $i$. When offering to sell energy to $\Acal$, this capacity remains unknown. Denote by $F$,  the joint cumulative distribution function (cdf) of 
$$\v{C} := (C_1, \ldots, C_N) \in [0, \bar{C}]^N.$$

DER supply being uncertain, $\Acal$ may not be able to supply $X$ in real-time that is promised in the DA market. Assume that $\Acal$ buys back the deficit $\left( X - \v{1}^\T \v{C} \right)^+$ at the real-time price $\LRT$, where we use the notation $z^+ := \max\{z, 0\}$ for a scalar $z$.  $\Acal$ then proceeds to allocate this penalty to the prosumers. We adopt the cost-sharing mechanism studied in \cite{costsharing} to design said penalties. Specifically, prosumer $i$ pays the penalty
\begin{align} 
\phi(x_i,\v{x}_{-i};\v{C})
&:= 
\LRT \(X - \v{1}^\T \v{C} \)^+ \frac{(x_i-C_i)^+}{\sum_{j=1}^N (x_j-C_j)^+}.
\label{eq:costsharing}
\end{align}
The penalty depends on her own offer $x_i$, the collective offers $\v{x}_{-i}$ of other prosumers and the realized supply $\v{C}$. \footnote{In practice, the price for the shortfall is not necessary the real-time market price $\LRT$. Our analysis remain unaffected if a different price for the penalty, call it  $\lambda$, is used.}

Such a penalty or cost sharing mechanism enjoys several desirable ``fairness'' axioms, thanks to the analysis in \cite{costsharing}. Each prosumer pays a fraction of the penalty that $\Acal$ pays for supply shortfall. A prosumer who is able to meet its promised supply does not pay any penalty. And, the penalty grows with the size of the shortfall. Finally, two prosumers with equal shortfalls face the same penalties. These properties are summarized mathematically in Figure \ref{fig:costsharing}. 
\begin{figure}
\centering
{\small\begin{mdframed}
\begin{itemize}[leftmargin=*]
\item Nonnegativity: $\phi(x_i,\v{x}_{-i};\v{C})\geq0.$
\item Budget balance: $\sum_{i=1}^N \phi(x_i,\v{x}_{-i};\v{C})= \LRT ( X - \v{1}^\T \v{C})^+.$
\item No exploitation: $x_i-C_i \leq 0 \implies \phi(x_i,\v{x}_{-i};\v{C})=0$.
\item Symmetry: $x_i-C_i = x_j-C_j \implies \phi(x_i,\v{x}_{-i};\v{C})= \phi(x_j,\v{x}_{-j};\v{C})$.
\item Monotonicity: $x_i-C_i \geq x_j-C_j \implies \phi(x_i,\v{x}_{-i};\v{C}) \geq \phi(x_j,\v{x}_{-j};\v{C})$.
\end{itemize}
\end{mdframed}}
\caption{Properties of the penalty sharing mechanism}
\label{fig:costsharing}
\end{figure}

Real-time transactions do not affect $\Acal$'s revenue. $\Acal$ therefore seeks a price $\rho$ that optimizes her profit from arbitrage in the day-ahead market and solves
 \begin{align}
\label{eq:aggregator.opt1}
\begin{alignedat}{5} 
\underset{\rho \geq 0}{\text{maximize}}& \ \ \pi_{\Acal}(\rho,\v{x}(\rho)):=(\lambda_{\text{DA}} - \rho) X(\rho).
	\end{alignedat}
\end{align}
In the above problem, we make explicit the dependency of $\v{x}$ and $X$ on $\rho$. 
Wholesale markets often only allow aggregations of a minimum size to participate, e.g., CAISO requires a minimum capacity of $0.5$ MW for a DER aggregator to participate \cite{CAISO}. Such restrictions can be included in \eqref{eq:aggregator.opt1}; we ignore them for ease of exposition. In fact, our previous analysis in \cite{CISS} showed that if one imposes a minimum capacity for the aggregator to participate, she might increase its offer price $\rho$ to prosumers, attracting them to offer a sufficiently large aggregate capacity.

Prosumer $i$ offers to sell $x_i$ amount of energy to $\Acal$. In doing so, he trades off between supplying  to $\Acal$ and consuming it locally. Denoting his utility of power consumption by $u$, prosumer $i$ solves
\begin{align}
\label{eq:prosumer.opt}
\begin{alignedat}{2}
 \underset{0\leq x_i \leq \bar{C}}{\text{maximize}} \ \  \pi_i(x_i,\v{x}_{-i},\rho) := \rho x_i + \E \left[ u\(d^0+C_i-x_i\) - \phi(x_i,\v{x}_{-i};\v{C}) \right],
\end{alignedat}
\end{align}
given $\Acal$'s offer price $\rho$. By selling $x_i$ to $\Acal$, she receives a compensation $\rho x_i$ at $\Acal$'s offer price $\rho$. 
Here, $d^0 \geq 0$ denotes the nominal energy consumption that prosumer $i$ purchases at a fixed retail rate either from a distribution utility or $\Acal$. We ignore the cost considerations of nominal demand as it does not affect our analysis.
Assume throughout that $u$ is nonnegative, concave, and increasing. Also, we let $d^0>\bar{C}$, i.e., DER supply is not large enough to cover the nominal demand.

Recall that prosumer $i$ offers $x_i$ before observing $C_i$ and therefore, the aggregate real-time supply from all $N$ prosumers can fall short of the promised supply. In that event, $\Acal$ faces a penalty that she allocates among the $N$ prosumers, $i$'s share being $\phi(x_i, \v{x}_{-i}; \v{C})$.
Prosumers do not believe their energy sales will affect real-time prices. Hence, the expected penalty for prosumer $i$ in \eqref{eq:prosumer.opt} is given by
\begin{align*}
& \E[\phi(x_i, \v{x}_{-i}; \v{C})]= \E[\LRT] \ \E\[ \(X - \v{1}^\T \v{C} \)^+ \frac{(x_i-C_i)^+}{\sum_{j=1}^N (x_j-C_j)^+} \]. 
\end{align*}
In the sequel, we abuse notation and write  $\LRT$ in place of $\E[\LRT]$ throughout. 


Given DA price $\LDA$ and expected real-time price $\LRT$, the prosumer-aggregator interaction can be summarized as a Stackelberg game $\Gfrak(\LDA, \LRT)$. $\Acal$ acts as a Stackelberg leader who decides price $\rho$. Prosumers in $\Nset$ follow by responding simultaneously with energy offers $\v{x}$. $\Acal$'s payoff is given by $\pi_\Acal$, while prosumer $i$'s payoff is $\pi_i$.
The pair $( \v{x}^*(\rho^*),\rho^*) $ constitutes a Stackelberg equilibrium of $\Gfrak(\lambda_{\text{DA}},\lambda_{\text{RT}})$, if 
$$\pi_i(x_i^*(\rho),\v{x}^*_{-i}(\rho), \rho) \geq \pi_i(x_i, \v{x}^*_{-i}(\rho),\rho)$$
for all $x_i \in [0, \bar{C}], \rho \geq 0, i \in \Nset$, and
$$\pi_{\Acal}\(\rho^*, \v{x}^*(\rho^*)\) \geq \pi_{\Acal}\(\rho,\v{x}^*(\rho)\)$$ 
for all $\rho \geq 0$.
We  establish in Theorem \ref{prop:NE} when such an equilibrium exists and is unique in the prosumer-aggregator game. The following assumption proves useful in the proof.

{\assumption[Random DER capacities]{
$F$ is smooth, fully supported on $[0,\bar{C}]^N$, and invariant under permutations.}
\label{assumption}
}

Smoothness and full support imply that DER capacities do not have any probability mass, and $F$ is strictly increasing in each argument over $[0,\bar{C}]^N$.
Invariance under permutations is natural for geographically co-located DERs, implying that prosumers' supplies are exchangeable. Given this assumption, $\Gfrak(\LDA, \LRT)$ becomes a symmetric game among the prosumers. We establish the existence of a Stackelberg equilibrium with symmetric reactions from prosumers. Later in this paper, we explicitly compute such equilibria and study their nature. 

{{\theorem[Existence and Uniqueness] 
Suppose Assumption \ref{assumption} holds. Then,  
a Stackelberg equilibrium 
$( \v{x}^*(\rho^*),\rho^*)$
always exists for $\Gfrak(\LDA, \LRT)$ with a unique symmetric response from prosumers, i.e., 
$x_i^*(\rho) = x^*(\rho)$
for each $\rho \geq 0$ and $i$ in $\Nset$. 
Furthermore, if 
\begin{equation} \frac{1}{2}\(\LDA - \rho\) \frac{\partial^2 X^*(\rho)}{\partial \rho^2}  <  \frac{\partial X^*(\rho)}{\partial \rho}
\label{eq:concavity} \end{equation}
 for all $\rho> 0$, then, the Stackelberg equilibrium with symmetric response from prosumers is unique, where $ X^*(\rho)=\v{1}^\T \v{x^*}(\rho)$.
\label{prop:NE}
}}

Our proof leverages a result from \cite{symmetric} that guarantees the existence of a symmetric Nash equilibrium of the game among prosumers, given $\Acal$'s price.\footnote{We remark that when prosumers are not homogenous and $F$ is not permutation invariant, Rosen's result in \cite{rosen} still guarantees the existence of a Nash equilibrium in the game among prosumers.} Exploiting the properties of the penalty sharing mechanism, we further establish that there exists a unique symmetric Nash equilibrium $\v{x}^*(\rho)$ that varies smoothly with $\rho$. $\Acal$ will never opt for a price higher than $\LDA$ and her profit varies smoothly in $\rho \in [0, \LDA]$. The smooth profit attains a maximum over that interval, leading to existence of a Stackelberg equilibrium. The relation in \eqref{eq:concavity} implies that $\Acal$'s profit becomes strictly concave and hence, the maximum and the equilibrium become unique. 

Next, we quantify the impact of DERs on wholesale market efficiency along the equilibrium path in the prosumer-aggregator game.

\section{Price of Aggregation}
\label{sec:efficiency}

Prosumers with available supply capacity can supplant conventional generation. Our goal is to characterize the impact of prosumer supply on the efficiency of the wholesale market. With a stylized wholesale market model, we compute the total energy procurement cost under two different models of prosumer participation. In the first model, aggregator $\Acal$ offers the aggregated procured supply capacity from individual prosumers to the wholesale market. The second model describes the ideal benchmark, where prosumers offer their supply capacity directly to the wholesale market. The comparison of the energy procurement costs in these two frameworks for prosumer participation leads to what we call the price of aggregation.


Consider a day-ahead wholesale market with dispatchable conventional generators and prosumers (participating directly or through an aggregator) competing to supply a point forecast of an inelastic demand $D$. Consider $G$ dispatchable generators, labelled $1, \ldots, G$. Let producer $j$ supply $Q_j$ amount of energy within its production capability set as $\[\underline{Q}_j,\overline{Q}_j\]$. Let its dispatch cost for producing $Q_j$ be given by $c_j(Q_j)$, where $c_j$ is a convex, nondecreasing and nonnegative function. We assume that $c_j$ truly reflects the production costs of the generator. In other words, we neglect possible market power of dispatchable power producers \cite{Baldick,Green, bose2014}, leaving a study of effects of strategic interactions of conventional generators and aggregated prosumer supply to future endeavors.
The economic dispatch problems in these two participation models that the system operator solves in day-ahead to clear the wholesale market are as follows.\footnote{We have ignored transmission network constraints in the wholesale market description for ease of exposition. Introducing these constraints causes no conceptual difficulty.
}

{When prosumers participate in the wholesale market through $\Acal$, the system operator clears the DA market by solving
\begin{eqnarray}
\label{eq.SO}
\Ccal^*_{\Acal} := & \underset{q_{\Acal},\v{Q}}{\text{minimum}} 
& \sum_{j=1}^G c_j(Q_j)  + \int_{0}^{q_{\Acal}} p_\Acal(y)dy, \nonumber \\
& \text{subject to}
&\underline{Q}_j \leq Q_j  \leq \overline{Q}_j, \ \ 0\leq q_{\Acal} \leq X, \\
&& \sum_{j=1}^G Q_j+q_{\Acal} = D.\nonumber
\end{eqnarray}
The inverse supply offer $p_\Acal(y)$ indicates the minimum price at which $\Acal$ is willing to sell $y$ amount of energy. The aggregated supply capacity for a wholesale market price $p_\Acal$ is given by $X[\rho^*(p_\Acal)]$, where $X=\v{1}^T\v{x}(\rho)$, and $(\v{x}(\rho^*),\rho^*)$ is the Stackelberg equilibrium of the game $\Gfrak(p_\Acal, \LRT)$.  By utilizing the equilibrium price path $\rho^*(p_\Acal)$ and taking the inverse of $X[\rho^*(p_\Acal)]$, one can compute the inverse supply offer $p_\Acal(X)$. Once \eqref{eq.SO} is solved, the day-ahead price  $\LDA$ is given by the optimal Lagrange multiplier of the supply-demand balance constraint.

{When prosumers directly participate in the wholesale market}, the system operator clears the DA market by solving
\begin{eqnarray}
\label{eq.SO2}
\Ccal_{\Pcal}^* :=&
 \underset{\v{q,Q}}{\text{minimum}} 
& \sum_{j=1}^G c_j(Q_j)  + \sum_{i=1}^N \int_{0}^{q_i} p_i(y_i)dy_i, \nonumber \\
& \text{subject to}
& \underline{Q}_j \leq Q_j  \leq \overline{Q}_j, \ \ 0\leq q_i \leq \bar{C}, \label{eq:SO.p}
\\
&& \sum_{j=1}^GQ_j+ \sum_{i=1}^Nq_i = D. \nonumber
\end{eqnarray}
Here, $p_i(y_i)$ is the inverse supply offer for each prosumer. Given a wholesale market price $p_i$ faced by prosumer $i$, his response $y_i(p_i)$ is computed by solving 
\begin{equation}\argmax_{y_i\in[0,\bar{C}]}  \Big\{ \ p_iy_i +  \E \[u\(d^0+C_i-y_i\)- \LRT(y_i-C_i)^+\] \Big\}. \label{eq:benchmark} \end{equation} Taking the inverse of $y_i(p_i)$, one can compute $p_i(y_i)$. The day-ahead price  $\LDA$ is the optimal Lagrange multiplier of the supply-demand balance constraint. Note that here, there is no cost sharing game among prosumers, as they are directly penalized for their corresponding shortfalls. Alternatively, one can take $y_i(p_i)$ to be the symmetric equilibrium response by prosumers $y(p)$, taking into account their cost shares, and then finding its inverse $p(y)$. In that case, the benchmark model is one in which DER supplies are concatenated via a purely social aggregator, who does not make profits from price arbitrage. Our analysis and insights in this paper are largely unaffected by such a variant of the benchmark model.



The supply capacity of a prosumer is typically too small for consideration in a wholesale market, and computing the dispatch and settlement for a large number of prosumers places an untenable computational burden on the system operator. To complicate matters, a transmission system operator typically does not have visibility into a distribution network. Hence, neither can they ensure that the DER dispatch will induce feasible flows in the distribution network, nor can they audit the actual supply. It is imperative that DER supply capacities are aggregated for participation in the wholesale market. The idealized direct prosumer participation model serves as a benchmark for the performance of any aggregation mechanism. 
Our other model for prosumer participation analyzes the case of a single profit-maximizing DER aggregator who chooses to represent the supplies from all prosumers in a system operator's footprint. In reality, such an entity will either be regulated or several aggregators will compete for prosumer representation, e.g., in \cite{mywork}. The loss in efficiency due to the strategic incentives of this single aggregator represents the maximum such loss the market will endure. Extending the model to incorporate competition for aggregation remains an interesting direction for future research.
Without DER participation, the market does not harness possible resources, and hence, is inefficient. However, the presence of $\Acal$ brings an efficiency loss due to the strategic incentives of $\Acal$, compared to the benchmark case in which prosumers participate directly in the wholesale market. We introduce the following metric of efficiency loss.
%
\begin{definition}
The Price of Aggregation ({\sf PoAg}) is given by $\frac{\Ccal^*_{\Acal}}{\Ccal^*_{\Pcal}}$, where $\Ccal^*_{\Acal}$ is computed using the Stackelberg equilibrium supply offer and $\Ccal^*_{\Pcal}$ is computed using \eqref{eq:SO.p}.
\end{definition}

$\text{\sf PoAg} \geq 1$ measures the efficiency loss of prosumer participation through an aggregator compared to direct prosumer participation. A larger {\sf PoAg} indicates a higher efficiency loss due to aggregation. We will now analyze {\sf PoAg} in various settings. To do so, one needs to construct the supply offers, that in turn, depends on the equilibrium in $\Gfrak(\LDA, \LRT)$. We construct these supply offers under two extreme cases for the stochasticity of DER supply.


\section{Studying the equilibrium of $\Gfrak(\LDA, \LRT)$.}\label{sec:game2}
While Theorem \ref{prop:NE} guarantees the existence of a Stackelberg equilibrium, it does not offer insights into the structure of said equilibrium. For general probability distributions on DER capacities in $\v{C}$, characterization of such equilibrium remains challenging. Here, we study two  settings for which such computation is easy--one where the capacities are completely dependent with identical components in $\v{C}$, and the other where the capacities are independent but identically distributed (iid). One expects the capacities in practice to follow a distribution that is somewhere  between these two extremes.

To motivate these two regimes, imagine that prosumers are supplying energy from rooftop solar. For geographically co-located prosumers, one expects high correlation among $C_i$'s. If two prosumers are geographically separated, independence among $C_i$'s might arise. A more realistic model is one with a collection of prosumer clusters that have highly correlated supply capacities within clusters, but independent between clusters. We relegate such considerations to future efforts and examine the two simpler settings here.

We study the iid case through the lens of mean-field (MF) games in the large prosumer limit $N \to \infty$. 
Such games have gained popularity following the seminal works in \cite{MF1,MF2} and are particularly useful to analyze interactions among a large number of players, where players respond to the population as a whole.
Taking  $N\rightarrow \infty$ in $\Gfrak(\LDA, \LRT)$ yields
 \begin{align}&\frac{\(\sum_{j=1}^N (x_j-C_j)\)^+ /N}{\sum_{j=1}^N  (x_j-C_j)^+/N} \rightarrow 
 \beta \ \text{almost surely},
 \label{eq:ratio}\end{align}
following the law of large numbers, where $\beta$ is a constant. Then, the penalty of shortfall for  prosumer $i$ becomes $\beta \LRT (x _i - C_i)^+$ and he maximizes 
 \begin{equation} 
 \E \[u\(d^0+C_i-x_i\)  + \rho x_i -  \beta \lambda_{\text{RT}} (x_i - C_i)^+ \],
 \label{MF_Problem}
 \end{equation}
given $\beta$. Further, $\beta$ is such that the solution of the above maximization satisfies \eqref{eq:ratio}.
Jensen's inequality on $f(z) = z^+$ yields $0 \leq \beta \leq 1$. Call this game $\Gfrak^\infty(\LDA, \LRT)$.

{{\theorem[Prosumer Offer Characterization]  Suppose  Assumption \ref{assumption} holds. If $\v{C}$ has identical components, i.e., $C_i = C$ for $i \in \Ncal$, then the prosumer offers for any $\rho \geq 0$ in $\Gfrak(\LDA, \LRT)$ satisfy
\begin{equation}
F(x^*(\rho)) = \frac{1}{\LRT}\( -\E\[u'(d^0+C-x^*(\rho))\]+\rho\).
\label{eq:NE2} 
\end{equation}
If $C_i$'s are iid, then in the limit $N\to \infty$, the prosumer offers $x^*(\rho)$ for any $\rho \geq 0$ in $\Gfrak^\infty(\LDA, \LRT)$ satisfy
\begin{align}
\beta F_i (x^*(\rho))   = \frac{1}{ \LRT}\( -\E\[ u'(d^0+C_i-x^* (\rho))\]+\rho\),
\label{eq:MF} 
\end{align}
where $\beta \in [0,1]$ is defined as
 \begin{equation} \beta =\frac{x^*(\rho) - \E[C_i] }{\E\[\(x^*(\rho) - C_i\)^+\]}. 
 \label{eq:beta}
 \end{equation}
\label{prop:MFNE}
}}

The expressions in \eqref{eq:NE2} and \eqref{eq:MF} are quite similar, implying that equilibrium offers of prosumers in $\Gfrak(\LDA,\LRT)$ behave similarly between the two extreme settings with completely dependent and independent DER capacities. These two settings provide the best and the worst-case {\sf PoAg}. The case with correlated but not fully dependent capacities in practice lies between these two extremes. 

As the capacities become independent, the solution to the mean-field game satisfies $x^*(\rho)\geq\E[C_i]$ owing to $\beta \geq 0$. That is, a prosumer may offer more than his mean anticipated capacity. This optimism stems from independence in supply capacities; a prosumer hopes that other prosumers will likely cover any shortfall on his part, leading to zero penalty. In contrast, prosumers with completely dependent capacities may choose $x^*(\rho)<\E[C_i]$ for small enough $\rho$. The room for optimism disappears as all prosumers face the same uncertainty. Equipped with these insights into the equilibrium in $\Gfrak(\LDA, \LRT)$, we  study the impact of DER aggregation on wholesale market efficiency through illustrative examples next.


\section{Illustrative Examples}\label{sec:examples}

We now study equilibrium offers of an example collection of prosumers, the same for an aggregator who participates in a wholesale market, and the resulting price of aggregation. Assume throughout that $N$ prosumers have linear utilities defined by $u(z)=\gamma z$ with $\gamma > 0$. \footnote{Accurately modeling the preferences $u(\cdot)$ of end-use customers in electricity consumption remains a challenging task (see \cite{feng} for a discussion) and it is beyond the scope of this paper. While we adopt linear utilities here, we remark that logarithmic utilities are commonly used in the economics literature \cite{shadow} for various commodities, and has recently found applications in the electricity market literature as well, e.g., in \cite{gao,sabita,mywork}. Using logarithmic utilities, along with deterministic DER supply, our previous analysis in \cite{CISS} reveals similar insights for price of aggregation. In this work, we adopt linear utilities for ease of exposition and offer insights into the impact of uncertainty. }


First, we examine the special case in which DER supply is deterministically $\bar{C}$, i.e., $C_i = \bar{C}$. Then, $\phi = 0$ and the payoffs of each prosumer is decoupled, given $\rho$. Each prosumer responds to $\rho$ by solving 
$$\underset{0 \leq x \leq \bar{C}}{\text{maximize}} \ \ u(d^0+\bar{C}-x)+\rho x.$$ 
If $\rho>\gamma$, each prosumer picks $x^*=\bar{C}$. As we show next, $\rho>\gamma$ is not sufficient for prosumers to sell their entire capacities when $\v{C}$'s are stochastic in nature. 

{\proposition{
If $C_i=C, \ i \in \Nset$ is uniformly distributed in $\[ \mu-\sqrt{3}\sigma, \mu+\sqrt{3}\sigma \]$, then the unique Stackelberg equilibrium $(\v{x}^*(\rho^*),\rho^*)$ of $\Gfrak(\LDA, \LRT)$ satisfies
\begin{align*}
x^*(\rho)&= \mu - \sqrt{3} \sigma +  \frac{2}{\LRT}(\rho-\gamma) \sqrt{3}\sigma, \ \text{for } \rho\geq \gamma,\\
\rho^*&=\frac{1}{2}(\LDA+\gamma)- \frac{\LRT(\mu-\sqrt{3}\sigma)}{4\sqrt{3}\sigma},
\end{align*}
where 
$$\sigma \in \[\frac{\LRT\mu}{2\sqrt{3}(\LDA-\gamma+\LRT/2)},\frac{\mu}{\sqrt{3}}\].$$
\label{prop:ex2}}}

In the presence of uncertainty, prosumers' response increases linearly in $\rho$, leading to the uniqueness of the Stackelberg equilibrium, similar to what we obtained in Theorem \ref{prop:NE}. 
Larger the average random capacity $\mu$, higher is the amount prosumers agree to sell. Higher the randomness, as captured by $\sigma$, lesser is the amount prosumers offer, owing to possible penalties they might face. Also, penalties are proportional to $\LRT$. Consequently, increasing $\LRT$ decreases their offers. Higher the possible penalty, either due to higher $\sigma$ or higher $\LRT$, aggregator needs to increase its price offer $\rho^*$ to attract DER supply. The bounds on $\sigma$ are the ones for which the equilibrium path is well-defined. 
\begin{figure}
\centering
\includegraphics[width=.75 \linewidth]{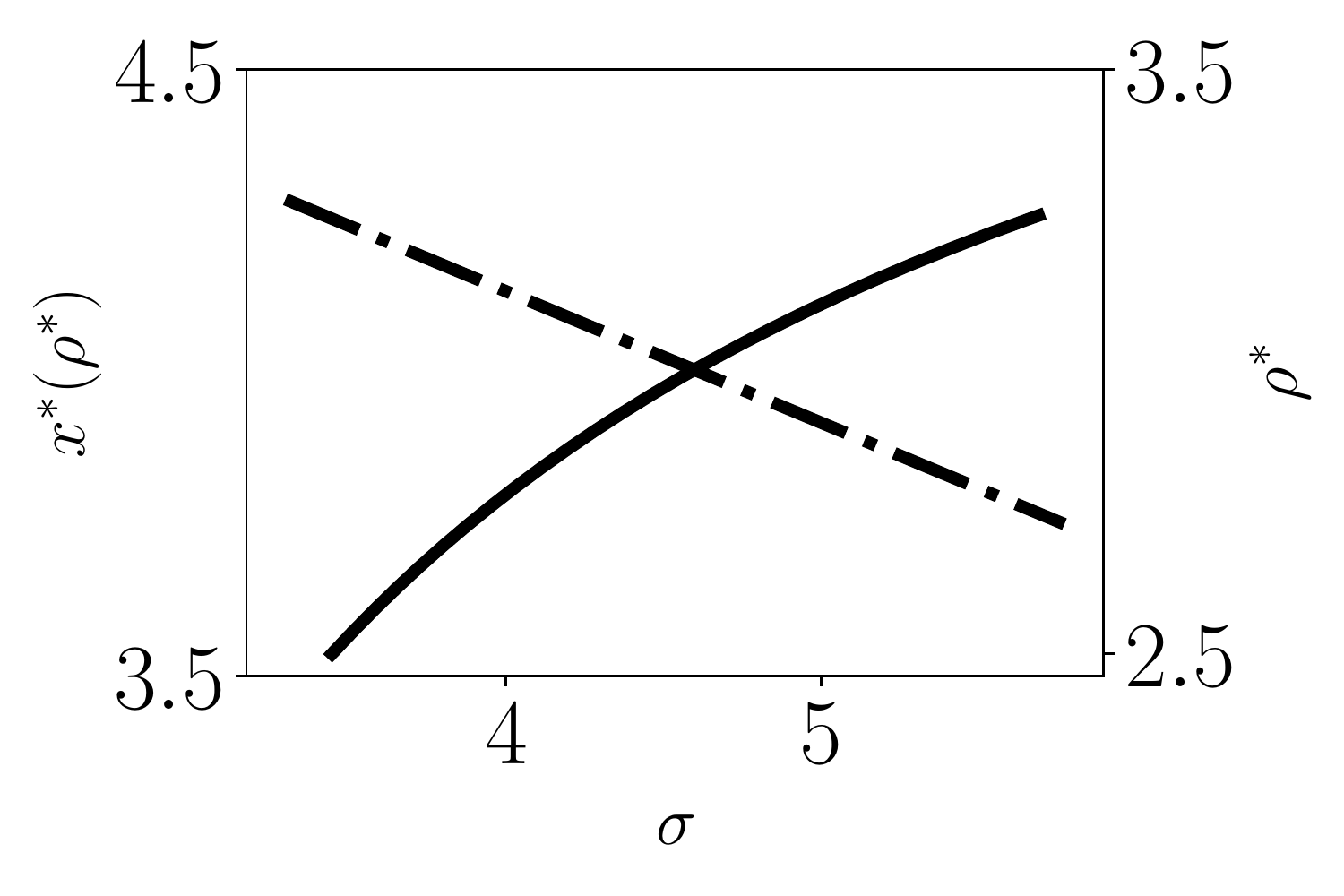}
\caption{Variation of equilibrium offer $x^*(\rho)$ and $\Acal$'s offer price $\rho^*$ with linear utilities. We choose $\gamma=2.5, \mu=10, \LDA=4, \LRT=4$.}
\label{illustration}
\end{figure}

Figure \ref{illustration} visualizes the effect of varying 
 $\sigma$. The plot corroborates our discussion above.
For $\sigma$ in the above range, we have $x^*(\rho)<\mu/2$. Contrast this result with iid capacities. Theorem \ref{prop:MFNE} reveals 
that the mean-field solution satisfies $x^*(\rho)\geq\mu$. In fact for this problem, 
we have $\beta = 0$ and $ x^*(\rho) = \mu$ for $\rho \geq \gamma$ and $\rho^*=\gamma$. 
The offer with iid capacities can be significantly higher than with identical capacities. 
%
%

\begin{figure}
\centering
  \includegraphics[width=\linewidth]{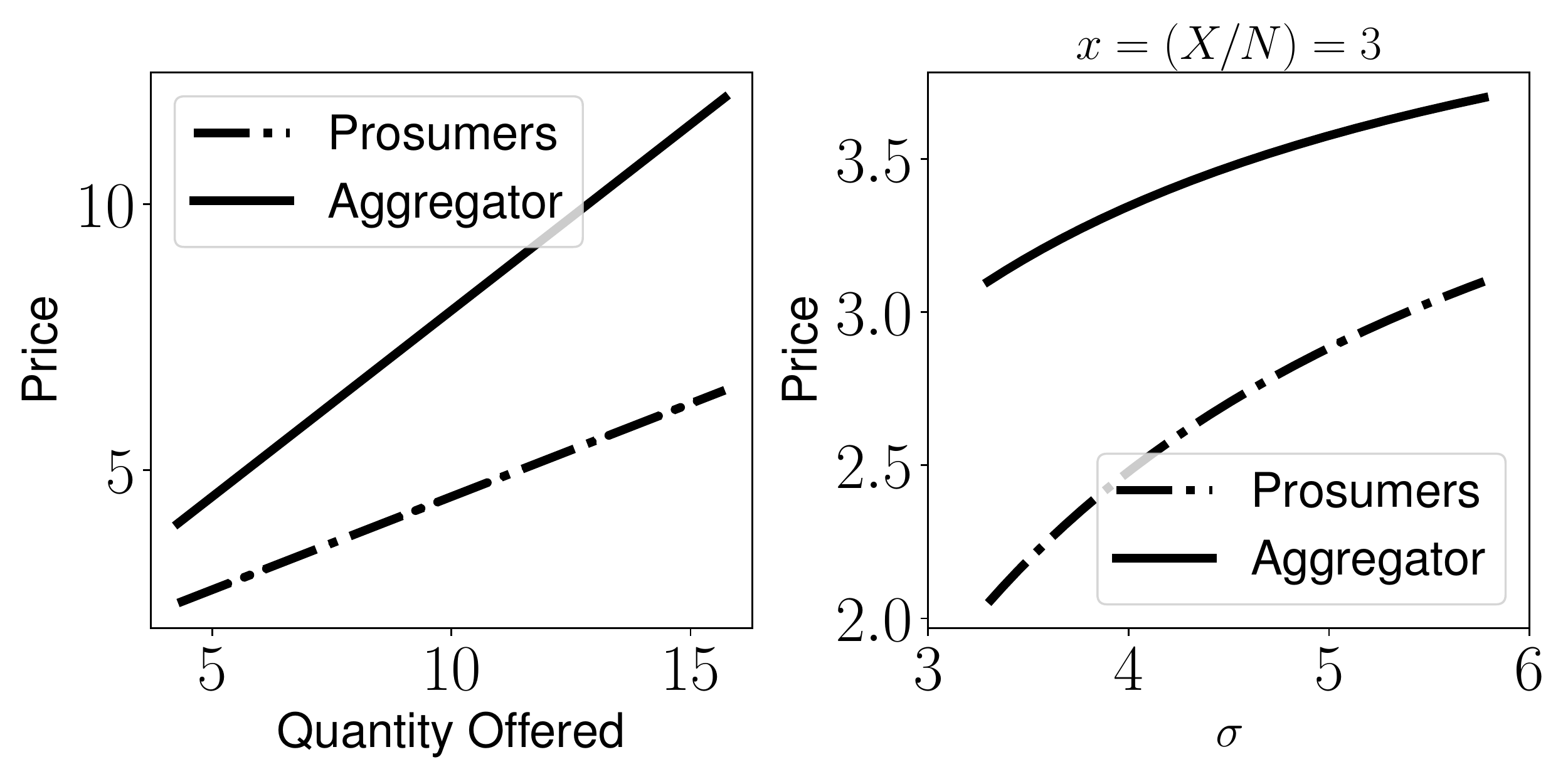}
  \caption{{\bf Left}: Equilibrium supply offers. Here, $\gamma=2.5, \mu=10,  \sigma=3.3,$ and $\LRT=4$. {\bf Right}: The influence of varying $\sigma$ on supply offers, for a fixed supply quantity.}
  \label{fig:supply_curves_stochastic}
\end{figure}
\subsection{Inverse Supply Functions}
The equilibrium of $\Gfrak(p_\Acal, \LRT)$ characterized in Proposition \ref{prop:ex2} implies the following aggregated supply as a function of $\Acal$'s price offer $p_\Acal$:
\begin{align}
X[\rho^*(p_\Acal)]=\frac{N}{2}\[\mu - \sqrt{3} \sigma +  \frac{2}{\LRT} \( {p_{\Acal}-\gamma} \)\sqrt{3}\sigma\] 
  \label{eq:optimal.x3}
\end{align}

The aggregator's inverse supply offer is then described by its inverse\begin{equation}
p_\Acal(X)=  \LRT  (2X/N -\mu+\sqrt{3}\sigma )/(2\sqrt{3}\sigma) + \gamma.\label{eq:inverse.x2}
\end{equation}

The best response provided by Proposition \ref{prop:ex2} analogously solves \eqref{eq:benchmark}, since capacities are completely dependent and prosumers are identical. The optimal supply offer by each prosumer is then
\begin{equation}
x(p)=\mu - \sqrt{3} \sigma +  \frac{2}{\LRT}(p-\gamma) \sqrt{3}\sigma.  \label{eq:optimal.x3}
\end{equation}

 Each prosumer's offer enters the objective function in problem \eqref{eq:SO.p} via its induced cost given by 
\begin{equation}
p(x)=\LRT(x- \mu+\sqrt{3}\sigma)/(2\sqrt{3}\sigma)  + \gamma . \label{eq:inverse.x}
\end{equation}

It is of importance to contrast how the aggregate supply capacity of the prosumers get offered in the wholesale market under two different prosumer participation models, using (\ref{eq:inverse.x2}) when $\Acal$ is present, and using (\ref{eq:inverse.x}) for the prosumers-only case.  Aggregator $\Acal$ offers the same supply capacity at a higher price than the collection of prosumers in aggregate. This price inflation is a consequence of the aggregator's aim to maximize her profits from arbitrage between the wholesale market prices and the prices she offers the prosumers (see Figure \ref{fig:supply_curves_stochastic}). To elaborate on the effects of uncertainty, we fix a quantity $X/N$, and then vary $\sigma$. We observe that as the variance increases, DER owners require higher prices that are also closer to $\Acal$'s offer price. This is a consequence of the risks of paying a penalty for the  shortfall. 


\subsection{The Price of Aggregation}
We now exploit the results in the previous subsection to analytically characterize the {\sf PoAg} in the next proposition.
 {\begin{proposition}
 {
In a wholesale market with $N$ DER suppliers, and one conventional generator with cost   
$$c(Q)=\kappa Q,\qquad \kappa>\gamma,\qquad Q\in[0,\infty],$$ the {\sf PoAg} is given by $\Ccal^*_{\Acal}/\Ccal^*_{\Pcal}$, 
where
\begin{align*} 
\Ccal_{\Acal}^*
&=\kappa D-\frac{Nq^*}{2} \( \kappa-\gamma+\frac{\LRT(\mu-\sqrt{3}\sigma)}{2\sqrt{3}\sigma} - \frac{q^*\LRT}{4\sqrt{3}\sigma} \), 
\\
\Ccal_{\Pcal}^*
&=\kappa D- {Nq^*}\( \kappa-\gamma+\frac{\LRT(\mu-\sqrt{3}\sigma)}{2\sqrt{3}\sigma} - \frac{q^*\LRT}{4\sqrt{3}\sigma}  \), 
\\
q^* &= \mu -  \sqrt{3}\sigma +\frac{2}{\LRT}\({\kappa-\gamma}\) \sqrt{3}\sigma,
\\
\sigma &\in \[\frac{\LRT\mu}{2\sqrt{3}(\LDA-\gamma+\LRT/2)},\frac{\mu}{\sqrt{3}}\].
\end{align*}
 Furthermore, in the absence of DER supply, the optimal procurement cost is $\kappa D.$
\label{prop:costs}
}\end{proposition}
}

One can readily observe by the above proposition that 
$$\kappa D>\Ccal_{\Acal}^*>\Ccal_{\Pcal}^*.$$ 

\begin{figure}
\centering
  \includegraphics[width=\linewidth]{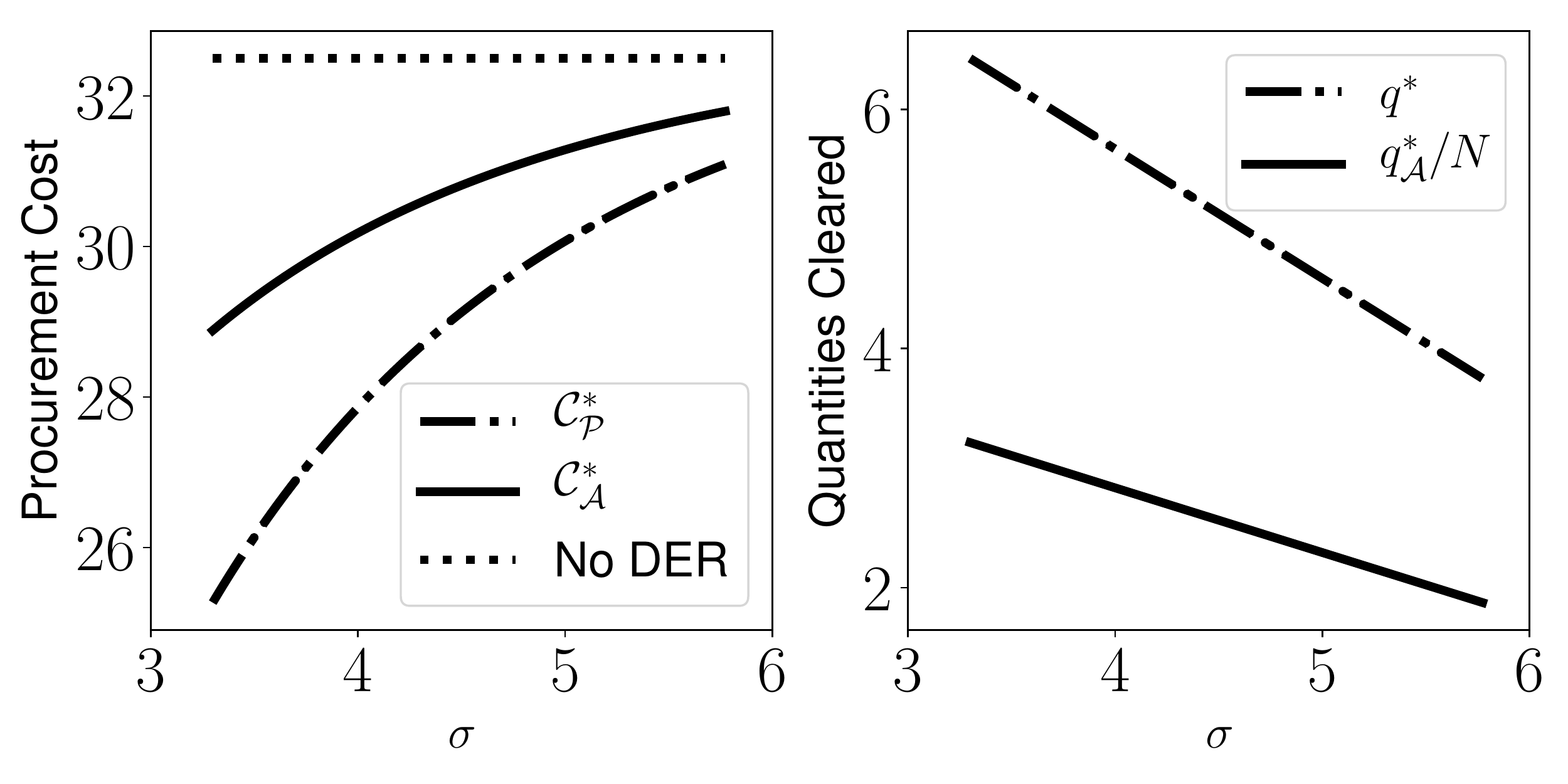}
  \caption{{\bf Left}: Procurement costs as $\sigma$ increases for $3$ cases: without DER supply (maximum cost), with DER supply and the existence of the aggregator (partial savings), and when prosumers offer their DER supply directly to the system operator (minimum cost). Increasing uncertainty makes DER participation unattractive to participants. {\bf Right}: DER quantities cleared as $\sigma$ increases.  More utilization of resources when DER owners directly participate. We use $\mu=10, \LRT=4, \gamma=2.5, \kappa=3.25, D/N=10$.}
  \label{fig:costs}
\end{figure}

\begin{figure}
\centering
  \includegraphics[width=\linewidth]{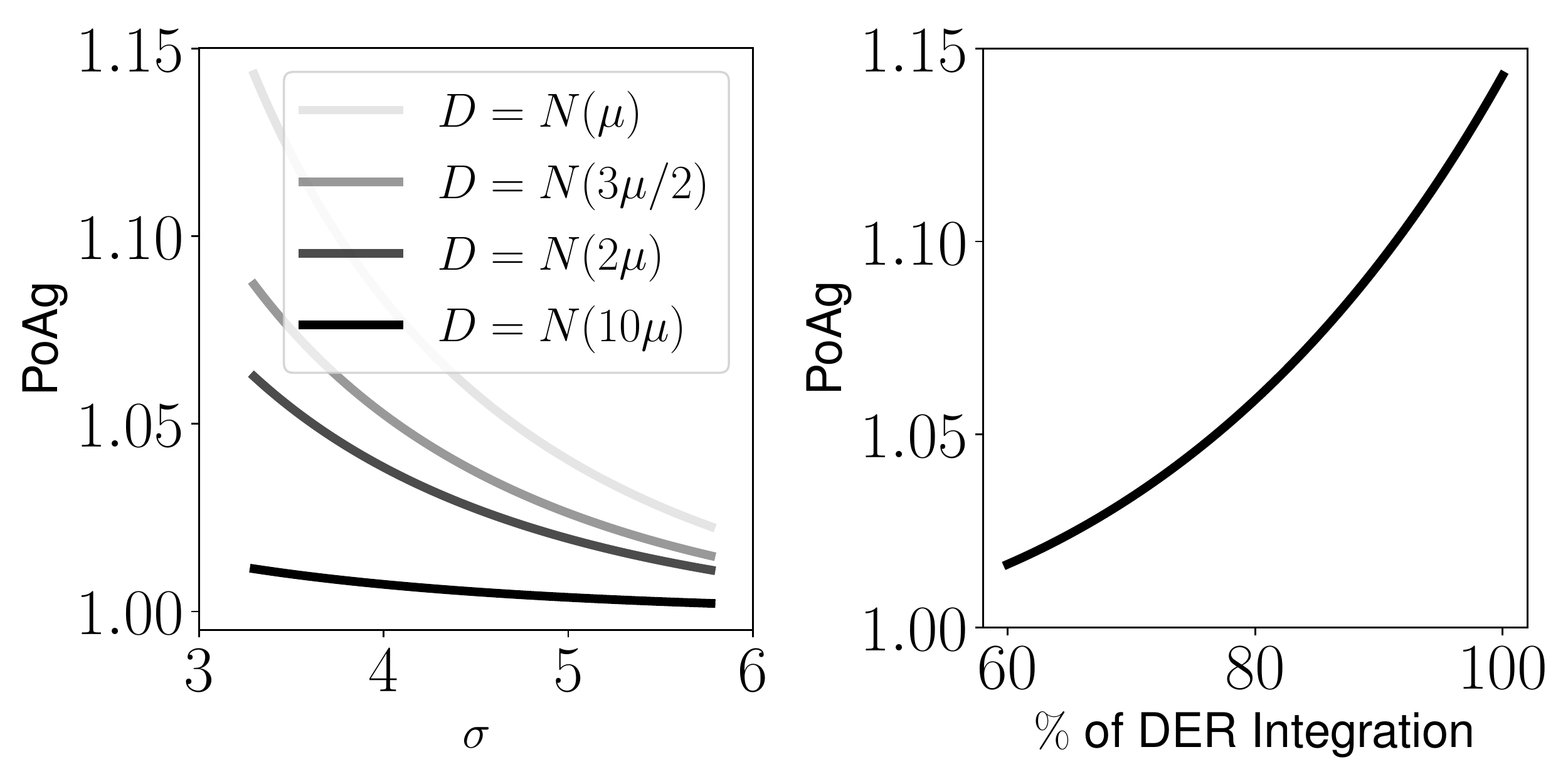}
  \caption{{\bf Left}: The price of aggregation as $\sigma$ increases. The {\sf PoAg} is monotonically decreasing with $\sigma$. {\sf PoAg} improves as DER contributes less towards the overall demand. {\bf Right}: The Price of Aggregation vs. $\%$ of DER supply. The {\sf PoAg} is monotonically increasing with DER integration. We use $\LRT=4, \gamma=2.5,$ and  $\kappa=3.25$. For the curve on the right, $D/N=10, \sigma=3.3,$ and $\mu$ varies from $6$ to $10$, reflecting percentage of integration.}
  \label{fig:PoAg}
\end{figure}

Figure \ref{fig:costs} plots the optimal procurement costs and  quantities cleared as $\sigma$ varies.  As one expects, procurement costs are minimized when prosumers offer their supply directly to the wholesale market, and savings diminish as uncertainty increases. Aggregation via profit-maximizing $\Acal$ strikes a balance between two extreme possibilities; no DER supply, and direct DER participation to the wholesale market.
This is also consistent with the supply curves in Figure \ref{fig:supply_curves_stochastic}. 
We note that intermediaries are inevitable, given the current wholesale market structures. 

As the uncertainty increases, {\sf PoAg} gets smaller as prosumers choose to sell less energy to the wholesale market, leaving $\Acal$ with smaller profits. As the DER integration increases, more prosumers sell energy, leading to increased profits for the aggregator. The worst-case {\sf PoAg} is attained at $100\%$ integration and $\sigma$ being near its lowest possible value at which the equilibrium path is well-defined.  Such a $\text{\sf PoAg}\approx1.15$ here, which implies that the cost with the aggregator is at most $15\%$ higher than the benchmark case. Both costs are still smaller than having no DER participation at all. Figure \ref{fig:PoAg} illustrates these tradeoffs.

\section{Conclusions and Future Directions} \label{sec:conclusion}
Our analysis points to debates surrounding the right design choice for incorporation of DERs in wholesale electricity markets. Should they be aggregated by third-party for-profit aggregators, perhaps where they vie to represent prosumers' supplies in the wholesale electricity market? Or, should a not-for-profit entity such as an independent distribution system operator be established to harness supply capacities of resources at the grid-edge? While the debates themselves are beyond the scope of this paper, we have provided a framework to quantify the benefits of different design choices.

In this paper, we considered and compared two different models of DER participation in wholesale electricity markets. In the first model, DERs directly offer their capacities in the wholesale market, while in the second one DERs participate in aggregate via a third-party for-profit aggregator. We modeled the strategic interactions between prosumers and the aggregator as a stochastic Stackelberg game. We characterized equilibria and explored two extreme cases: DER capacities are completely dependent or independent and identically distributed. At the equilibrium, we quantified the effects of aggregation through a metric we called Price of Aggregation ({\sf PoAg}). 

There are several directions for future work. We assumed the DER aggregator to be a price-taker in the wholesale market. However, under high penetration of DERs, it may be possible that aggregators can influence the wholesale market price and engage in strategic bidding. 
Furthermore, wholesale electricity markets typically have multiple settlements for energy procured in each hour. Analysis of DER participation with and without an aggregator with stochastic supply and multi-settlement wholesale market structure is another important direction for future work.

\bibliographystyle{IEEEtran}
\bibliography{references}
\section{Appendix}
\subsection{Proof of Theorem \ref{prop:NE}}
The proof proceeds in four steps. First, we show that a symmetric Nash equilibrium exists for the game among prosumers. Then, we show that this symmetric equilibrium is unique. The analysis then allows us to conclude that a Stackelberg equilibrium exists. We finally identify  sufficient conditions under which a Stackelberg equilibrium with symmetric response from prosumers is unique.

\emph{Step 1: Existence of symmetric Nash equilibrium}. For each prosumer $i$, we show that $\pi_i(x_i,\v{x}_{-i},\rho)$ is concave in $x_i\in[0,\bar{C}]$ for each $\v{x}_{-i}$, and it is jointly continuous in $(x_i,\v{x}_{-i})$ to conclude that a pure-strategy Nash equilibrium exists in the concave Nash game among the prosumers. Furthermore, under Assumption \ref{assumption}, the game is symmetric, and we infer the existence of a symmetric Nash equilibrium from \cite{symmetric}. Towards that goal, it suffices to show that the expected cost share $\E[\phi(x_i,\v{x}_{-i};\v{C})]$ is both smooth and convex in $x_i$ for each $\v{x}_{-i}$ and continuous in $\v{x}_{-i}$, in view of the expression in \eqref{eq:prosumer.opt}. 
We begin by defining
$$s_{-i}:=\sum_{j\neq i} (x_j-C_j), \quad S_{-i}:=\sum_{j\neq i} (x_j-C_j)^+.$$ 
Then, the derivative $\[\frac{\E[\phi(x_i,\v{x}_{-i};\v{C})]}{\lambda_{\text{RT}}}\]^{'}$ w.r.t. $x_i$ is 
\begin{align}
&\Bigg[\E \[\frac{\(s_{-i}+x_i-C_i\)^+}{S_{-i}+(x_i-C_i)^+} (x_i-C_i)^+ \] \Bigg]^{'} \nonumber\\ 
&\qquad= \E\[\(\frac{(s_{-i}+x_i-C_i)(x_i-C_i)}{(S_{-i}+x_i-C_i)} \mathbbm{1}_{\{C_i\leq\bar{C}_i\}}\)^{'}\] \nonumber \\
&\qquad=\E\[\( 1 + \frac{S_{-i}(s_{-i}-S_{-i})}{(S_{-i}+x _i-C_i)^2}\)\mathbbm{1}_{\{C_i\leq\bar{C}_i\}}\], \label{eq:der_phi}
\end{align}
where $\bar{C}_i:=\min\{x_i,x_i+s_{-i}\}$, and we have used completion of squares to get the second equality. 

The quotient of two continuous functions is continuous, provided that the denominator is non-zero. Thus, it follows that $\[\E[\phi(x_i,\v{x}_{-i};\v{C})]\]^{'}$ is continuous in $x_i$ since $F$ is smooth. Hence, $\E[\phi(x_i,\v{x}_{-i};\v{C})]$ is smooth. Analogously, one can verify that $\E[\phi(x_i,\v{x}_{-i};\v{C})]$ is continuous in $\v{x}_{-i}$.  Using \eqref{eq:der_phi}, we then have
 \begin{equation*}
 \[\frac{\E[\phi(x_i,\v{x}_{-i};\v{C})]}{\lambda_{\text{RT}}}\]^{''}  = 
\E\[\frac{2S_{-i}(S_{-i}-s_{-i})\mathbbm{1}_{\{C_i<\bar{C}_i\}}}{(S_{-i}+x _i-C_i)^2}\] \label{eq:convex}
\end{equation*}
is nonnegative since $S_{-i} \geq s_{-i}$ and $S_{-i}+x _i-C_i > 0$. Hence, $\E[\phi((x_i,\v{x}_{-i},\v{C})]$ is convex in $x_i$ for each $\v{x}_{-i}$. Following \cite[Theorem 3]{symmetric}, a symmetric equilibrium exists.

 \emph{Step 2: Uniqueness of symmetric Nash equilibrium.}
Karush-Kuhn-Tucker optimality conditions and \eqref{eq:der_phi} imply that the best response $x_i^*$, given $\v{x}_{-i}$ and $\rho$, satisfies
\begin{align}
&\frac{\E\[ -u'(d^0+C_i-x_i^*)\]+\rho + \underline{\eta} -\bar{\eta}}{\lambda_{\text{RT}}} =\E\[\( 1 + \frac{S_{-i}(s_{-i}-S_{-i})}{(S_{-i}+x^* _i-C_i)^2}\)\mathbbm{1}_{\{C_i\leq\bar{C}_i\}}\], \label{eq:KKT_1}
\end{align}
where $\underline{\eta}$ is the nonnegative Lagrange multiplier associated with the constraint $x_i^*\geq0$ and $\bar{\eta}$ is the same associated with $x_i^*\leq \bar{C}$. 
Notice that $S_{-i}\geq s_{-i}$ and hence, we have two possible events for which the right hand side of \eqref{eq:KKT_1} can be nonzero -- when $S_{-i}=s_{-i}$ and when $S_{-i}>s_{-i}$. Denote these events by $\hat{\Xcal}_i$ and $\Xcal_i$, where
\begin{align*}
\hat{\Xcal}_i:&= \{ \v{C} \ \mid \ 0\leq C_j\leq x_j, \ j\in\Nset, \ s_{-i}=S_{-i} \}, \\
\Xcal_i  :&=\{  \v{C} \ \mid \ 0 \leq C_i\leq\bar{C}_i, s_{-i}<S_{-i} \}.
\end{align*}
Then, we have 
\begin{align*}
&\E\[\( 1 + \frac{S_{-i}(s_{-i}-S_{-i})}{(S_{-i}+x^* _i-C_i)^2}\)\mathbbm{1}_{\{\hat{\Xcal}_i\}}\]=F(x^*_i,\v{x}_{-i}).\end{align*}
%
\begin{align*}
\E\[\( 1 + \frac{S_{-i}(s_{-i}-S_{-i})}{(S_{-i}+x^* _i-C_i)^2}\)\mathbbm{1}_{\{\Xcal_i\}}\]=:h_i(x^*_i,\v{x}_{-i}). 
\end{align*}
%
%
A Nash equilibrium is attained at the intersection of best responses, and hence, by the first-order conditions, $\v{x}^*(\rho)$ at any Nash equilibrium solves
\begin{align}
&\frac{\E\[ -u'(d^0+C_i-x_i^*)\]+\rho+ \underline{\eta} -\bar{\eta}}{\lambda_{\text{RT}}}  =F(x_i^*,\v{x}^*_{-i}) + h_i(x^*_i,\v{x}^*_{-i}) \label{eq:NashKKT}
\end{align}
for  each $i\in\Nset$.  Symmetric equilibria of the Nash game among prosumers satisfy the above equation with $x^*_i(\rho)=x^*(\rho)$ for each prosumer $i$. Substituting $x^*$ for $x_{j}^*$ for $j \in \Nset$ in said relation, we get
\begin{equation}\frac{\E\[ -u'(d^0+C_i-x^*)\]+\rho+ \underline{\eta} -\bar{\eta}}{\lambda_{\text{RT}}} = F(x^*) + h(x^*),\label{eq:NashKKT.s}\end{equation}
where we slightly abuse notation and write
$$F(x^*) := F(x^*, \ldots, x^*), \quad h(x^*) = h_i(x^*, \ldots, x^*).$$
The above notation for $h$ is well-defined, owing to the permutation invariance of the distribution on $\v{C}$.
%
We now show that for each $\rho$, there is a single $x^*$ that solves \eqref{eq:NashKKT.s}. To do this, we utilize the following observations: 
\begin{itemize}[leftmargin=*]
\item {\bf Observation 1:} $ 0\leq F(x^*) + h(x^*) \leq 1$. This follows from properties of $\E[\phi]$.
\item {\bf Observation 2:} $h(x)$ is increasing in $x$. To see this, notice that
$$\frac{\partial s_{-i}}{\partial x} = N-1, \qquad  \frac{\partial S_{-i}}{\partial x} =: M < N-1, $$
where $N$ is the number of prosumers. Then, we have
{\small
\begin{align*}
h'(x)&=\E\[\(\frac{S_{-i}(s_{-i}-S_{-i})}{(S_{-i}+x-C_i)^2}\)' \mathbbm{1}_{\{ \Xcal_i\}}\]\\
&=\E\[\frac{M(s_{-i}-S_{-i})+S_{-i}(N-1-M)}{(S_{-i}+x-C_i)}\mathbbm{1}_{\{ \Xcal_i\}}\] - \E\[\frac{2S_{-i}(M+1)(s_{-i}-S_{-i})}{(S_{-i}+x-C_i)^2}\mathbbm{1}_{\{ \Xcal_i\}}\]\\
&=\E\Bigg[\frac{\mathbbm{1}_{\{ \Xcal_i\}}}{(S_{-i}+x-C_i)^2} \times \Big( 2(M+1)((S_{-i})^2-s_{-i}S_{-i}))+ (Ms_{-i}+S_{-i}(N-1-2M)) (S_{-i}+x-C_i) \Big) \Bigg]\\
&\geq0,\end{align*}
}
where the inequality follows from $N-1>M$, $x-C_i\geq0$, and $S_{-i}>s_{-i}$ for each $ \Xcal_i$. 
\item {\bf Observation 3:} $\underline{\eta}$ and $\overline{\eta}$ cannot be both strictly positive. This is a straightforward consequence of complementary slackness conditions for the prosumer's problem. 

\end{itemize}

Equipped with these observations, define 
\begin{align*}\rho_{\min}&:=- \E\[ -u'(d^0+C_i-x_i)\]\Big\vert_{x_i=0}, \\
\rho_{\max}&:=\lambda_{\text{RT}}-  \E\[ -u'(d^0+C_i-x_i)\]\Big\vert_{x_i=\bar{C}}
\end{align*}
and split the analysis of possible symmetric equilibria into five cases, depending on $\rho$.
\begin{itemize}[leftmargin=*]

\item When $\rho< \rho_{\min}$: From the properties of $u$, it follows that $\E\[ -u'(d^0+C_i-x^*)\]+\rho < 0$ for all $x^*\in[0,\bar{C}]$. From \eqref{eq:NashKKT.s}, we infer that $\underline{\eta}>0$, that in turn implies  $x^*=0$ and $\bar{\eta}=0$ by complementary slackness condition.

\item When $\rho= \rho_{\min}$: For all $x^*\in[0,\bar{C}]$, by the above observations, we must have $\bar{\eta}=0$, and hence,
\begin{align}&\frac{\E\[ -u'(d^0+C_i-x^*)\]- \E\[ -u'(d^0+C_i)\]+ \underline{\eta} }{\lambda_{\text{RT}}} = F(x^*) + h(x^*). 
\label{eq:rho.rho.min}\end{align}
If $\underline{\eta}>0$, then $x^*=0$, and we have $$\frac{\underline{\eta}}{\LRT} -  (F(0) + h(0)) = \frac{\underline{\eta}}{\LRT} = 0,$$ which is a contradiction. Plugging $\underline{\eta}=0$ in \eqref{eq:rho.rho.min} yields
$x^*=0$ as the only solution in view of the  the properties of $F$, $h$, and $u$.

\item When $\rho_{\min} < \rho <\rho_{\max}$: 
Define
\begin{equation} g(x) := \frac{\E\[ -u'(d^0+C_i-x)\]+\rho}{\lambda_{\text{RT}}} - F(x) - h(x). \label{eq:NashKKT2}\end{equation}
Then, the properties of $F$, $h$, and $u$ imply that $g$ is strictly decreasing with $g(0)>0$ and $g(\bar{C})<0$. Hence, it has a unique zero-crossing in $(0,\bar{C})$ that identifies $x^*$.

\item When $\rho=\rho_{\max}$: 
For all $x^*\in[0,\bar{C}]$, we must have $\underline{\eta}=0$. Also, we have
\begin{align*}&\frac{\E\[ -u'(d^0+C_i-x^*)\]+ \lambda_{\text{RT}} -  \E\[ -u'(d^0+C_i-\bar{C})\] -\bar{\eta} }{\lambda_{\text{RT}}} = F(x^*) + h(x^*).\end{align*}
If $\bar{\eta}>0$, then $x^*=\bar{C}$ and $F(\bar{C}) + h(\bar{C}) = 1$. Then, we get
$( \lambda_{\text{RT}} -\bar{\eta} )/\lambda_{\text{RT}}   =  1$
and $\bar{\eta}=0$, yielding a contradiction. Plugging $\bar{\eta}=0$ in the above equation yields 
$x^*=\bar{C}$ as the only solution.

\item When $\rho>\rho_{\max}$: For each $x^*\in[0,\bar{C}]$, we have 
$$\frac{\E\[ -u'(d^0+C_i-x^*)\]+\rho}{\lambda_{\text{RT}}}>1,$$
which implies that $\bar{\eta}>0$, and hence $x^*=\bar{C}$ and $\underline{\eta}=0$.

\end{itemize}
This completes the proof of unique symmetric best response from the prosumers for each $\rho$. Combining this observation with step 1 yields the existence of a unique symmetric Nash equilibrium.

\emph{Step 3: Existence of Stackelberg equilibrium. }
\indent From the above analysis, we infer that $x^*(\rho)$ is both unique and continuous in $\rho$, and so is  $X^*(\rho) = N {x}^*(\rho)$. Thus, $\pi_{\Acal}(\rho,\v{x}^*(\rho))$ is  continuous, and the set over which $\rho$ can take values is compact. 
Hence, existence of a solution to \eqref{eq:aggregator.opt1} is guaranteed by  Weierstrass Theorem \cite{Luenberger}, implying that a Stackelberg equilibrium always exists. 

\emph{Step 4: Uniqueness of Stackelberg equilibrium with symmetric prosumer response. }
Notice that
$$\frac{\partial^2  \pi_{\Acal}(\rho,\v{x}^*(\rho)) }{\partial \rho^2}=-2 \frac{\partial X^*(\rho)}{\partial \rho} + (\lambda_{\text{DA}}-\rho)  \frac{\partial^2 X^*(\rho)}{\partial \rho^2}.$$ 
When the inequality holds, $\pi_{\Acal}(\rho,\v{x}^*(\rho))$ is strictly concave in $\rho$, leading to the uniqueness of the symmetric Stackelberg equilibrium strategy $\rho^*$, completing the proof.

\subsection{Proof of Theorem \ref{prop:MFNE}}

(a) At the symmetric equilibrium for completely dependent capacities, it follows that
$$ \Xcal_i=\emptyset,  \ i\in\Nset \implies  h(x^*)=0.$$ Via the analysis of $\eqref{eq:NashKKT.s}$ in the proof of Theorem \ref{prop:NE}, it follows that it is optimal for $\Acal$ to pick $\rho_{\min}\leq\rho\leq\rho_{\max}$ because choosing $\rho<\rho_{\min}$ or $\rho>\rho_{\max}$ does not increase $\Acal$'s profit. It also follows via the same analysis that the symmetric equilibrium is given by solving $$ F(x^*(\rho)) = \frac{1}{\LRT}\(\E\[-u'(d^0+C-x^*(\rho))\]+\rho\),$$
proving the claim.

(b) By the first-order condition of the mean-field problem in \eqref{MF_Problem}, each player solves 
\begin{align}
&\E\[ -u'(d^0+C_i-x^*)\]+\rho  =\lambda_{\text{RT}}\beta\[\int_{0}^{x^*} (x^*-C_i) f(C_i)dC_i\]^{'}=\beta F_i(x^*), \label{eq:x_MF}
\end{align}
where the last equality follows from the Leibniz integral rule. Note that $\beta \geq 0$ by definition, and by the law of large numbers, $\beta$ solves  
$$\beta=\frac{\(\E\[x^*- C_i\]\)^+}{\E\[\(x^* - C_i\)^+\]},$$
yielding the desired result. 
\subsection{Proof of Proposition \ref{prop:ex2}} 

Analogous to the derivation of  \eqref{eq:NashKKT.s} in the proof of Theorem \ref{prop:NE}, linear utilities and completely dependent capacities with $\rho\geq\gamma$ and  the Karush-Kuhn-Tucker optimality conditions imply 
$$ F\(x^*(\rho)\)= \frac{1}{\lambda_{\text{RT}}}(\rho-\gamma).$$ 
 When $C$ is uniformly distributed over $\[\mu-\sqrt{3}\sigma, \mu+\sqrt{3}\sigma\]$, it follows that 
$$ x^*(\rho)=\mu+ \(\frac{2}{\lambda_{\text{RT}}} \(\rho-\gamma\) -1\)\sqrt{3}\sigma, \ \ \rho\geq\gamma.$$

Plugging the above into the aggregator's problem, $\Acal$'s payoff is 
$$\pi_{\Acal}(\rho,\v{x}^*(\rho))= N(\lambda_{\text{DA}}-\rho)\[\mu+ \(\frac{2}{\lambda_{\text{RT}}} \(\rho-\gamma\) -1\)\sqrt{3}\sigma\],$$
which is quadratic and strictly concave in $\rho$. Taking the derivative and solving $\frac{\partial \pi_{\Acal}(\rho,\v{x}^*(\rho))}{\partial\rho}=0$ leads to the statement of the proposition, where $\sigma$ is picked such that the equilibrium path is well-defined. 
nd{cases}
%
%
\subsection{Proof of Proposition \ref{prop:costs}}
Since $\Acal$ is an inframarginal supplier, it follows that $\lambda_{\text{DA}}=\kappa.$ With the aggregator, using (\ref{eq:inverse.x2}), and by the supply-demand balance constraint, $Q=D-q_{\Acal}$, it follows that $\Ccal^*_{\Acal}$ is given by
$$ \min_{q_{\Acal}\in[0,N\bar{C}]} \[\kappa D+\frac{\lambda_{\text{RT}}q_{\Acal}^2}{2N\sqrt{3}\sigma} +\(\gamma -\kappa -\frac{(\mu-\sqrt{3}\sigma)\lambda_{\text{RT}}}{2\sqrt{3}\sigma}\)q_{\Acal}\].$$
Solving the above yields
$$q^*_{\Acal}=\frac{N}{2} \[\mu-  \sqrt{3}\sigma+\frac{2}{\LRT}\({\kappa-\gamma}\)  \sqrt{3}\sigma\].$$

 Similarly, without the aggregator, $\Ccal^*_{\mathcal{P}}$ is given by
$$ \min_{q\in[0,\bar{C}]} \[\kappa D+ \frac{N\lambda_{\text{RT}}q^2}{4\sqrt{3}\sigma} +\(\gamma -\kappa -\frac{(\mu-\sqrt{3}\sigma)\lambda_{\text{RT}}}{2\sqrt{3}\sigma}\)Nq\],$$
which is solved at 
$$ q^*= \mu-  \sqrt{3}\sigma+\frac{2}{\LRT}\({\kappa-\gamma}\)  \sqrt{3}\sigma = 2q^*_{\Acal}/N.$$
With no DER supply at all, $Q^*=D$, and the  cost is $\kappa D$. The statement of the proposition follows from simple re-arrangements, and taking into account the range for $\sigma$ such that the equilibrium path is well-defined. 

\end{document}